\documentclass[a4paper]{article}

                            \usepackage{amsmath,amssymb,graphics}
                             \usepackage{fancyhdr,epic}
                             \usepackage[latin1]{inputenc}
                             \newcommand{\nc}{\newcommand}
                             \newcommand{\rnc}{\renewcommand}

                             \nc{\nin}{\noindent}
                             \nc{\ssk}{\smallskip}
                             \nc{\msk}{\medskip}
                             \nc{\bsk}{\ssk\msk}
                             \nc{\sni}{\ssk\nin}
                             \nc{\mni}{\msk\nin}
                             \nc{\bni}{\bsk\nin}
                             \nc{\psk}{\msk}
                             \nc{\q}{\quad}
                             \nc{\qq}{\qquad}
                             \nc{\lb}{\linebreak}
                             \nc{\up}{\vspace{-\smallskipamount}}

                             \nc{\ul}[1]{\underline{#1}}
                             \nc{\ol}[1]{\overline{#1}}
                             \nc{\N}{\mathbb{N}}
                             \nc{\Z}{\mathbb{Z}}
                             \nc{\Q}{\mathbb{Q}}
                             \nc{\R}{\mathbb{R}}
                             \nc{\C}{\mathbb{C}}
                             \nc{\K}{\mathbb{K}}
                             \rnc{\P}{\mathbb P}
                             \nc{\F}{\mathbb{F}}

\nc{\quo}[2]{\raisebox{3pt}{$#1$}/\raisebox{-3pt}{$#2$}}

                             \nc{\x}{\cdot}
                             \nc{\ph}{\varphi}
                             \nc{\eps}{\varepsilon}
                             \nc{\tht}{\vartheta}
                             \nc{\eq}{\,=\,}
                             \nc{\bel}{\!\in\!}
                             \nc{\p}{^{\,\prime}}
                             \nc{\varempty}{{\rm

O}\hspace{-0.65em}\raisebox{0.18ex}{/}\hspace{0.22em}}
                             \nc{\inc}{\!\subseteq\!}
                             \nc{\us}{\underset}
                             \nc{\os}{\overset}
                             \nc{\cvd}{\hfill$\Box$}
                             \nc{\ve}[2]{{#1}_1,\ldots,{#1}_{#2}}
                             \nc{\ft}[2]{#1,\ldots,#2}

                             \nc{\dph}{\vphantom{\Bigg|}}
                             \nc{\dsum}[2]{\us{#1}{\os{#2}{\sum}}}
                             \nc{\dprod}[2]{\us{#1}{\os{#2}{\prod}}}
                             \nc{\dint}[2]{{\displaystyle\int_{#1}^{#2}}}
                             \nc{\dlim}[1]{\us{#1}{\lim}}
                             \nc{\dmax}[1]{\us{#1}{\max}}
                             \nc{\dmin}[1]{\us{#1}{\min}}
                             \nc{\comb}[2]{\left(\!\!\begin{array}{c}#1\\
                             #2\end{array}\!\!\right)}

                             \nc{\ct}[2]{\put(#1){\makebox(0,0){#2}}}
                             \nc{\cst}[2]{\put(#1){\makebox(0,0){\small #2}}}

                             \nc{\spdg}{sen\-za per\-di\-ta di
ge\-ne\-ra\-li\-tà}

                             \nc{\defi}{\nin{\bf Definizione.}\q}
                             \nc{\propr}{\nin{\bf Proprietà:}\q}
                             \nc{\teor}{\nin{\bf Teorema.}\q}
                             \rnc{\emph}[1]{{\sl #1}}
                             \nc{\prmq}{\us{{v\in M(\Q)}}{\prod}}
                             \nc{\prmk}{\us{{v\in M(K)}}{\prod}}
                             \nc{\prmqo}{\us{{v\in
M(\Q)\backslash\{0\}}}{\prod}}
                             \nc{\prmqi}{\us{{v\in
                             M(\Q)\backslash\{\infty\}}}{\prod}}
                             \nc{\prmqp}{\us{{v\in

M(\Q)\backslash\{\infty,p_0,\ldots,p_n\}}}{\prod}}
                             \nc{\prvp}{\us{{v\in
\{p_0,\ldots,p_n\}}}{\prod}}
                             \nc{\conv}{\os{n\to\infty}
                                 {-\!\!\!-\!\!\!-\!\!\!\!\!\longrightarrow}}
                             \nc{\convm}{\os{m\to\infty}
                                 {-\!\!\!-\!\!\!-\!\!\!\!\!\longrightarrow}}
                             \nc{\convy}{\os{y\to\infty}
                                 {-\!\!\!-\!\!\!-\!\!\!\!\!\longrightarrow}}
                             \nc{\convvax}{\os{|x|\to\infty}
                                 {-\!\!\!-\!\!\!-\!\!\!\!\!\longrightarrow}}
                             \nc{\nconvvax}{\os{|x|\to\infty}

{-\!\!\!-\!\!\!-\!\!\!\!\!\not\longrightarrow}}
                             \nc{\convvay}{\os{|y|\to\infty}
                                 {-\!\!\!-\!\!\!-\!\!\!\!\!\longrightarrow}}
                             \nc{\cD}{{\cal D}}
                             \nc{\cU}{{\cal U}}
                             \nc{\cP}{{\cal P}}
                             \nc{\ula}{\ul{\alpha}}
                             \nc{\ut}[1]{\ul{\widetilde{#1}}}
                             \nc{\oQ}{\ol{\Q}}
                             \nc{\MCD}{\hbox{\footnotesize MCD}}
                             \nc{\hK}{\widehat{K}}
                             \nc{\hv}{\widehat{v}}
                             \nc{\cu}{c_1(d)}
                             \nc{\cd}{c_2(d)}
                             \nc{\ep}{e^{2\pi i\vartheta}}
                             \nc{\ai}{\{a_i\}}
                             \nc{\bi}{\{b_i\}}
                             \nc{\aip}{\{a_i'\}}
                             \nc{\lii}{\dlim{i\to+\infty}}
                             \nc{\lni}{\dlim{n\to+\infty}}
                             \rnc{\L}{\Lambda}
                             \nc{\ub}{\ul{b}}
                             \nc{\vol}{{\rm vol}}
                             \nc{\dfracc}[2]{\dfrac{\raisebox{-2pt}{$#1$}}
                             {\raisebox{1pt}{$#2$}}}
                             \nc{\A}{{\rm Area}}
                             \nc{\ox}{\ol{x}}
                             \nc{\oy}{\ol{y}}
                             \nc{\sD}{\sqrt{D}}
                             \nc{\LV}{Liouville}
                             \nc{\xy}{\dfracc{x}{y}}
                             \nc{\yx}{\dfracc{y}{x}}
                             \nc{\pq}{\dfracc{p}{q}}
                             \nc{\ab}{\dfracc{a}{b}}
                             \nc{\vaaxy}{\Big|\alpha-\xy\Big|}
                             \nc{\he}{\widehat{\eps}}
                             \nc{\Bigva}[1]{\Big|#1\Big|}
                             \nc{\Bigvau}[1]{\hbox{${\Big|#1\Big|}_{11}$}}
                             \nc{\biggva}[1]{\bigg|#1\bigg|}
                             \nc{\biggvau}[1]{\hbox{${\bigg|#1\bigg|}_{11}$}}
                             \nc{\vau}[1]{\hbox{${|#1|}_{11}$}}
                             \nc{\rtd}{\sqrt[3]{2}}
                             \nc{\rtum}{\sqrt[3]{\dfracc12}}
                             \nc{\Lvi}{L_{v,i}}
                             \nc{\ux}{\ul{x}}
                             \nc{\ksa}{K\Sigma_A}
                             \nc{\ksap}{K\Sigma_A^+}
                             \rnc{\l}{\!\ll\!}
                             \nc{\ddsum}[2]{\us{#2}{\dsum{#1}{}}}
                             \nc{\ddprod}[2]{\us{#2}{\dprod{#1}{}}}
                             \nc{\bcn}{b_1\,c_1^n}
                             \nc{\abc}{congettura $a\,b\,c$}
                             \rnc{\Im}{{\rm Im}}
                             \nc{\twodots}{\hspace{-0,1em}.\:.\:}
                             \nc{\sqx}{\sqrt{x}}
                             \nc{\od}{\ol{d}}
                             \newtheorem{theorem}{Theorem}[section]
\newtheorem{lemma}[theorem]{Lemma}

\newtheorem{corollary}[theorem]{Corollary}
\newtheorem{Main Theorem}[theorem]{Main Theorem}
                             
\begin{document}
                             \title {On the period of the continued fraction for values of the square root of power
sums}

\author{Amedeo {\sc SCREMIN\footnote{\nin The author was supported by Istituto Nazionale di Alta
Matematica ``Francesco Severi'', grant for abroad Ph.D}}}


                             \date {\today}
                 \setlength{\baselineskip}{20.736pt}
                             \setlength{\parindent}{2.25em}
                             \maketitle

\begin{abstract}

\nin The present paper proves that if for a power sum $\alpha$
over $\Z$ the length of the period of the continued fraction for \
$\sqrt{\alpha(n)}$ \ is constant for infinitely many even (resp.
odd) $n$, then \ $\sqrt{\alpha(n)}$ \ admits a functional
continued fraction expansion for all even (resp. odd) $n,$ except
finitely many; in particular, for such $n,$ the partial quotients
can be expressed by power sums of the same kind.

\end{abstract}

\bni

\section{Introduction}

\nin It is well known that the continued fraction for rational
numbers is finite and that for the square root of a positive
integer \ $a$ \ which is not a square is periodic of the form \
$[a_o;\overline{a_1, \ldots, a_{R-1},2a_0}]$ \ (here with \
$\overline{a_1, \ldots, a_{R-1},2a_0}$ \ we denote the periodic
part), where \ $R\geq 1$ \ is the length of the period. About $R,$
\ we know that the bound \ $R\ll \sqrt{a}\log a$ \ holds (see
\cite {Hua} and \cite {L}).

\mni A power sum \ $\alpha$ \ is a function on \ $\N$ \ of the
form
\begin{equation}\label{intro1}
\alpha(n) = b_1 c_1^n + b_2 c_2^n + \ldots + b_h c_h^n,
\end{equation}
\nin where the roots \ $c_i$ \ are integers and the coefficients \
$b_i$ \ are in $\Q$ or in $\Z.$ We know from Corollary 1 in
\cite{C-Z} that, apart from the case when $\alpha$ is the square
of a power sum of the same kind, \ $\sqrt{\alpha(n)}$ is a
quadratic irrational for all \ $n \in \N,$ \ except finitely many.
\ This means that the continued fraction expansion for
$\sqrt{\alpha(n)}$ \ is periodic for $n$ large, raising the
problem weather the length of the period is bounded or not for
$n\longrightarrow +\infty,$ which will be considered in this
paper. Some partial results on such problem have been recently
obtained by Bugeaud and Luca (see \cite{B-L}).

\nin On a similar problem, but considering a non constant
polynomial \ $f$ \ with rational coefficients instead of the power
sum \ $\alpha,$ \ remarkable results were obtained by Schinzel in
\cite {Sc1} and \cite {Sc2}. He provided conditions on $f$ under
which the length of the period of the continued fraction for \
$\sqrt{f(n)}$ \ tends to infinity as $n\longrightarrow +\infty$.

\mni In the present paper we shall first prove that if a power sum
$\alpha$ with rational coefficients cannot be approximated "too
well" by the square of a power sum of the same kind, then the
length of the period of the continued fraction for
$\sqrt{\alpha(n)}$ tends to infinity as $n\longrightarrow +\infty$
\ (Corollary \ref{cor1}).

\nin Then we shall consider power sums with integral coefficients,
and show that for any fixed $r\in\{0,1\}$, \ if the length of the
period of the continued fraction for $\sqrt{\alpha(2m+r)}$ is
constant for all $m$ in an infinite set, then for every $m \in
\N,$ except finitely many exceptions, the partial quotients of the
continued fraction for $\sqrt{\alpha(2m+r)}$ can be identically
expressed by power sums of the same kind (Main Theorem
\ref{cor2}).

\nin The results above shall be deduced from some lower bounds for
the quantities \ $|\sqrt{\alpha(n)} - \dfrac{p}{q}\big|$ \
(Corollary \ref{teo1}) \ and \ $\Big|\dfrac{\sqrt{\alpha(n)} +
\beta(n)}{\gamma(n)} - \dfrac{p}{q}\Big|$ \ (Theorem \ref{teo2})
respectively, \ where \ $\alpha, \beta, \gamma$ \ are power sums
and \ $p,q$ \ are integers, which we shall obtain using Schmidt's
Subspace Theorem in a way similar to that of Corvaja and Zannier
in \cite{C-Z} and \cite{C-Z2}.

\nin Theorem \ref{teo2} and Corollary \ref{teo1} (taking
$\alpha=0$ and $q=1$ respectively), are the analogue of the
Theorem in \cite{C-Z2} and of Theorem 3 in \cite{C-Z}.

\nin The results contained in this paper give an answer to some
questions raised in the Final Remark (b) in \cite{C-Z2}, where it
is predicted that "under suitable assumptions on the power sum
$\alpha$ with rational roots and coefficients, the length of the
period of the continued fraction for $\sqrt{\alpha(n)}$ tends to
infinity with $n$".

\mni

\section{Notation}

\nin In the present paper we will denote by \
                             $\Sigma$ \ the ring of functions on $\N$,
                             called {\sl power sums}, of the
                             form
\begin{equation}\label{not1} \alpha(n) = b_1 c_1^n + b_2 c_2^n +
\ldots +b_hc_h^n,
\end{equation}                              \
\nin where the distinct {\sl roots} \ $c_i\neq 0$ \ are
                             in $\Z$, and the
                             {\sl coefficients} \ $b_i \in \Q^\star$. For rings $A,B\subseteq \C,$ let
                             $A\Sigma_B$ denote the ring of power
                             sums with coefficients in $A$ and
                             roots in $B.$

                             \nin If $B \subseteq \R$, it is usually enough to deal with power sums with only positive roots.
                             Working in
                             this domain causes no loss of generality: the
                             assumption of positivity
                             of the roots may usually be achieved by writing \ $2n+r$ \ instead
                             of \
                             $n$, and considering the cases of \ $r=0,1$
                             \ separately.

                             \sni If $\alpha\in\oQ\Sigma_\Q$, we
                             set
                             $l(\alpha):=\max\{c_1,\ldots,c_h\}$.
                             In the same way we define the
                             function $l$ for a power sum defined on the sets of even
                             or odd numbers.
                             It is immediate to check that \
                             $l(\alpha\beta)=l(\alpha)l(\beta)$, \
                             $l(\alpha+\beta)\leq
                             \max\{l(\alpha),l(\beta)\}$ \ and
                             that \
                             $l(\alpha)^n\gg |\alpha(n)|\gg
                             l(\alpha)^n$.

\mni {\bf NOTE} \ In the statements of our results and in the
proofs we will always omit the condition for the existence of
$\sqrt{\alpha(n)}\in \R$, i.e. that $\alpha(n)\geq 0$ for $n$
large.

\mni

\section{Statements}

\nin The following Theorem \ref{teo2} states that for power sums
$\alpha, \beta, \gamma \in \Sigma,$ if
$\dfrac{\sqrt{\alpha}+\beta}{\gamma}$ cannot be well approximated
on the subsequence of even (or odd) numbers by a power sum in
$\Sigma,$ then $\dfrac{\sqrt{\alpha(n)}+\beta(n)}{\gamma(n)}$
cannot be well approximated by rationals with exponentially
bounded denominators, except for a finite number of even (odd)
$n$. \ This Diophantine approximation result will be obtained
using Schmidt's Subspace Theorem in a way similar to that of
Corvaja and Zannier in \cite{C-Z} and \cite{C-Z2}. Theorem
\ref{teo2} is the main tool we will use to prove the Corollaries
and the Main Theorem.

\begin{theorem}\label{teo2} \q Let \ $\alpha, \beta, \gamma
\in \Sigma$, \ $\gamma$ not identically zero, and let \ $\eps>0$
and $r\in\{0,1\}$ be fixed.

\nin Suppose that there does not exist a power sum \
$\eta\in\Sigma$ \ such that
$$\Big|\dfrac{\sqrt{\alpha(2m+r)} +
\beta(2m+r)}{\gamma(2m+r)}-\eta(m)\Big|\ll e^{(-2m+r)\epsilon}.$$
\nin Then there exist $k=k(\alpha, \beta, \gamma)>2$ and
$Q=Q(\epsilon)>1$ with the following properties. For all but
finitely many naturals \ $n \equiv r \mod 2$ \ and for integers
$p,q$,\ $0<q<Q^{2m+r},$ we have
\begin{equation}\label{eqteo3}
\Big|\dfrac{\sqrt{\alpha(n)} + \beta(n)}{\gamma(n)} -
\frac{p}{q}\Big|\geq \frac{1}{q^k} e^{-\epsilon n}.
\end{equation}

\end{theorem}

\bni {\bf Remark 1} \q Taking $\alpha=0$ in Theorem \ref{teo2}, we
obtain again the result of the Theorem in \cite{C-Z2}.

\bni Corollary \ref{teo1} is a simplified version of Theorem
\ref{teo2}. It states that if a power sum $\alpha \in \Sigma$
cannot be well approximated on the subsequences of even and odd
numbers by the square of a power sum from the same ring, then
$\sqrt{\alpha(n)}$ cannot be well approximated by rationals with
exponentially bounded denominators, except for a finite number of
$n$. It will be used to prove Corollary \ref{cor1}.

\begin{corollary}\label{teo1} \q Let \ $\alpha\in\Sigma$, and let \ $\eps>0$ \ be
fixed. Assume that for every \ $r\in\{0,1\}$ \ and for all \
$\xi\in\Sigma$,
$$l(\alpha - \xi^2)\geq l(\alpha)^{1/2}$$
on the sequence $n=2m+r.$

\nin Then there exist \ $k=k(\alpha)> 2$ \ and \ $Q=Q(\epsilon)>1$
\ with the following properties. For all but finitely many
$n\in\N$ and for all integers \ $p,q, \ 0<q<Q^n,$ we have
\begin{equation}\label{eqteo1}
\big|\sqrt{\alpha(n)} - \frac{p}{q}\big|\geq \frac{1}{q^k}
e^{-\epsilon n}.
\end{equation}
\end{corollary}

\bni {\bf Remark 2} \q Taking \ $q=1,$ \ we can see that Corollary
\ref{teo1} is a generalization of Theorem 3 in \cite{C-Z}.

\sni {\bf Remark 3} \q In concrete cases, it is easy to verify
whether the assumption of Corollary \ref{teo1} holds or not. In
fact, it is enough to prove that for every $r\in \{0,1\}$ and for
all $\eta \in \Sigma$, \ in the power sum \
$\alpha(2m+r)-\eta(m)^2$ \ there cannot be cancellations of all
the coefficients of the roots greater than the square root of the
dominating root of $\alpha$ (resp., there exists $\eta$ such that
we have all that cancellations). By a similar way it is possible
to verify if the assumption of Theorem \ref{teo2} holds or not.

\bni The following Corollary \ref{cor1} states that if a power sum
$\alpha \in \Sigma$ cannot be well approximated by the square of a
power sum of the same kind, then the length of the period of the
continued fraction for $\sqrt{\alpha(n)}$ tends to infinity as
$n\longrightarrow +\infty.$ This result was already obtained with
a similar proof by Bugeaud and Luca in \cite[Theorem 2.1]{B-L}.

\begin{corollary}\label{cor1} \q Let $\alpha\in \Sigma$ be as
in the Corollary \ref{teo1}.

\nin Then the length of the period of the continued fraction for
$\sqrt{\alpha(n)}$ tends to infinity as $n\rightarrow +\infty.$

\end{corollary}

\bni The Main Theorem \ref{cor2} follows again from Theorem
\ref{teo2}, and states that if the length of the period of the
continued fraction for the square root of a power sum is constant
for infinitely many even (resp. odd) $n,$ \ then the partial
quotients of the continued fraction can be expressed by power sums
for all even (resp. odd) $n,$ \ except finitely many.

\begin{Main Theorem}\label{cor2} \q Let \ $\alpha\in
\Z\Sigma_\Z$, \ and let \ $r\in\{0,1\}$ \ be fixed.

\nin Suppose that there exists an infinite set $A\subseteq \N$ and
a constant $R\geq 0$ such that for \ $m\in A$ \ the length of the
period of the conti\-nued fraction expansion for
$\sqrt{\alpha(2m+r)}$ \ is \ $R.$

\nin Then there exist \ $\beta_0, \ldots, \beta_R \in \Z\Sigma_\Z$
\ such that for every \ $m\in \N,$ \ apart from finitely many
exceptions, we have the continued fraction expansion
\begin{equation}\label{eqcor2}
\sqrt{\alpha(2m+r)}=[\beta_0(m);\overline{\beta_1(m),\ldots,\beta_R(m)}].
\end{equation}

\end{Main Theorem}

\bni {\bf Remark 4} \q The result of Corollary \ref{cor1},
together with the Main Theorem \ref{cor2}, gives an answer to the
question raised in the Final Remark (b) in \cite{C-Z2}.

\mni

\section{Auxiliary results}

\nin For the reader's convenience we state here a version of
Schmidt's Subspace Theorem due to H.P. Schlickewei; we have
borrowed it from \cite[Theorem 1E, p. 178]{S2} (a complete proof
requires also \cite{S1}). It will be our main tool to prove
Theorem \ref{teo2}.

\begin{theorem}\label{sub} \q Let $S$ be a finite set of absolute values of $\Q$, including the
infinite one and normalized in the usual way (i.e. $|p|_v=p^{-1}$
if $v|p$). Extend each $v\in S$ to $\oQ$ in some way. For $v\in S$
let $L_{1,v},\ldots,L_{n,v}$ be $n$ linearly independent linear
forms in $n$ variables with algebraic coefficients and let $\delta
>0.$

\nin Then the solutions $\underline{x}:=(x_1,\ldots,x_n)\in \Z^n$
to the inequality
$$\dprod {v\in S}{} \ \dprod{i=1}{n}
|L_{i,v}(\underline{x})|_v<\dmax {1\leq i\leq n}|x_i|^{-\delta}$$
\nin are contained in finitely many proper subspaces of $\Q^n.$

\end{theorem}

\bni The following Lemma \ref{lem1} is a result by Evertse (in a
more general case); a proof by Corvaja and Zannier can be found in
\cite[Lemma 2]{C-Z}.

\begin{lemma}\label{lem1} \q Let \ $\xi\in \Sigma_\Q$ and let $D$ be the minimal
positive integer such that $D^n\xi\in\Sigma$.

\nin Then, for every \ $\eps>0,$ \ there are only finitely many \
$n\in\N$ \ such that the denominator of \ $\xi(n)$ \ is smaller
than \ $D^ne^{-n\eps}.$

\end{lemma}

\mni

\section{Proofs}

\nin We start with the following very simple

\begin{lemma}\label{lem3} \q Let \ $\alpha, \beta, \gamma \in \Sigma$, \
$\gamma$ not identically zero, and let \ $t$ \ be any positive
real number. Then for every \ $r\in \{0,1\}$ \ there exists \
$\eta_r \in \oQ\Sigma_{\Q}$ \ such that
$$\Big|\dfrac{\sqrt{\alpha(2m+r)} + \beta(2m+r)} {\gamma(2m+r)}
-\eta_r(2m+r)\Big| \ll t^{2m}.$$
\nin Such $\eta_r$ can be
effectively computed in terms of \ $r$, \ $\alpha, \beta, \gamma $
and $t$.
\end{lemma}
\mni {\bf Proof of Lemma \ref{lem3}.} \q Let \
$\alpha(n)=\dsum{j=1}{h}\,b_j\,c_j^n$, \
                             with \ $c_j \in
                             \Z,$ $c_j\neq0$ and \ $b_j \in \Q^*$ \ $\forall \ j=1,\ldots,h.$

\nin We can suppose $c_1>c_2>\ldots>c_h>0$.

\nin For a real determination (resp. real positive) of $b_1^{1/2}$
(resp. $c_1^{1/2}$), fixed for the rest of the proof, we have
\begin{equation}\label{eq0lemma1}
\alpha(n)^{1/2}=(b_1c_1^n)^{1/2}\Big(1+\dsum{j=2}{h}\dfrac{b_j}{b_1}\big(\dfrac{c_j}{c_1}\big)^n\Big)^{1/2}\eq
(b_1 c_1^n)^{1/2}(1+\sigma(n))^{1/2},
\end{equation}
\nin with $\sigma(n) \in \Sigma_\Q$, and $\sigma(n) \eq
O((c_2/c_1)^n)$.

                             \nin Expanding the function \ $x \mapsto
                             (1+x)^{1/2}$ \ in Taylor
                             series, we have
                             \begin{equation}\label{eq1lemma1}
                             (1+\sigma(n))^{1/2} \eq 1+
                             \dsum{j=1} {H} a_j \
                             \sigma(n)^j
                             +O(|\sigma(n)|^{H+1}),
                             \end{equation}
                             \nin where $H>0$ is an integer that can be chosen
                             later and $a_j$,\ $j=1,\ldots,H,$ \ are the Taylor
                             coefficients \ $
\big({1/2 \atop j}\big)$ \ of the
                             function \ $x \mapsto (1+x)^{1/2}$.

\nin For every $r\in \{0,1\},$ \ substituting (\ref{eq1lemma1}) in
(\ref{eq0lemma1}) we obtain
\begin{equation}\label{eq2lemma1}
\alpha(2m+r)^{1/2}=b_1^{1/2}c_1^{r/2}c_1^m
                             \Big(1+ \dsum{j=1} {H} a_j \sigma(2m+r)^j\Big)
                             +O\Big(\Big(\dfrac{c_2}{c_1}\Big)^{2m(H+1)}c_1^m\Big).
                             \end{equation}
\nin Let
\begin{equation}\label{eq4lemma1}
\beta(n)=\dsum{j=1}{k}d_je_j^n \in \Sigma ,
\end{equation}
\nin with \ $e_j \in
                             \Z,$ $e_j\neq0$ and \ $d_j \in \Q^*$ \ $\forall \ j=1,\ldots,h.$

\nin We can suppose $e_1>e_2>\ldots>e_k>0$.

\nin Fix $H$ such that \
$\Big(\dfrac{c_2}{c_1}\Big)^{(H+1)}c_1^{1/2}<e_1.$

\sni Let $\gamma(n)=\dsum{j=1}{l}f_jg_j^n\in\Sigma,$ with \ $g_j
\in \Z,$ $g_j\neq0$ and \ $f_j \in \Q^*$ \ $\forall \
j=1,\ldots,h.$

\nin We can suppose $g_1>g_2>\ldots>g_k>0$.

\nin Using the same method as in the proof of Theorem 1 in
\cite{C-Z}, we can write
\begin{equation}\label{eq3lemma1}
\gamma(n)^{-1}=f_1^{-1}g_1^{-n} \dsum{j=0}{s}
\phi(n)^j+O((g_2/g_1)^{n(s+1)}g_1^{-n}),
\end{equation}
\nin where $\phi(n):
=-\dsum{i=2}{l}\frac{f_i}{f_1}\big(\frac{g_i}{g_1}\big)^n \in
\Sigma_\Q$, \ $\phi(n)=O(g_2/g_1)^n,$ and $s>0$ is an integer that
can be chosen later.

\sni Thus, by equations (\ref{eq2lemma1}), \ (\ref{eq4lemma1}),
(\ref{eq3lemma1}), by the choice of $H$ and the definition of
$\phi$, we obtain

\msk $\dfrac{\sqrt{\alpha(2m+r)} + \beta(2m+r)}
{\gamma(2m+r)}=f_1^{-1}g_1^{-r}g_1^{-2m}\Big(\dsum{i=0}{s}\phi(2m+r)^{i}\Big)\cdot$

\sni
$\cdot\Big(b_1^{1/2}c_1^{r/2}c_1^m\big(1+\dsum{i=1}{H}a_i\sigma(2m+r)^{i}\big)+\dsum{i=1}{k}d_ie_i^{2m+r}\Big)
+ O\Big(\big(g_2/g_1\big)^{2m(s+1)}g_1^{-2m}e_1^{2m}\Big).$

\mni Fix $s$ such that $\big(g_2/g_1\big)^{(s+1)}g_1^{-1}e_1<t$ \
and put, for $r=0,1,$

\msk
$\eta_r(2m+r):=f_1^{-1}g_1^{-r}g_1^{-2m}\Big(\dsum{i=0}{s}\phi(2m+r)^{i}\Big)\cdot$

\q \q \q $\cdot
\Big(b_1^{1/2}c_1^{r/2}c_1^m\big(1+\dsum{i=1}{H}b_i\sigma(2m+r)^{i}\big)+
\dsum{i=1}{k}d_ie_i^{2m+r}\Big).$

\mni By definition \ $\eta_r\in \oQ\Sigma_\Q$ \ for every \
$i=0,1.$

\nin Thus for every \ $r\in \{0,1\}$ \ we have effectively
constructed a power sum \ $\eta_r(n)\in\oQ\Sigma_{\Q}$ \ such that

\msk $\Big|\dfrac{\sqrt{\alpha(2m+r)} + \beta(2m+r)}
{\gamma(2m+r)} -\eta_r(2m+r) \Big|\ll t^{2m}$,

\sni completing the proof.

\cvd

\nin {\bf Remark 5} \q Let us notice that in \ $\eta_r$ \ the root
with largest absolute value is \ $g_1^{-2}\cdot \max\{e_1^2,
c_1\}$ and that the other roots appearing are rational with
denominator powers of $c_1$ and $g_1.$ The denominators of each of
such roots are divided by $g_1^2.$

\bni {\bf Proof of Theorem \ref{teo2}.} \q Let $\eta_r,$ for
$r\in\{0,1\}$ fixed, be as in Lemma \ref{lem3}, with $t=1/9$.

\nin We can write (recall Remark 6)
$$\eta_r(2m+r)=b_{1,r}^{1/2}d_1^m
\big(g^{-2m}+b_2d_2^{2m+r}+\ldots+b_hd_h^{2m+r}\big),$$
\nin for
some $b_{1,r}, b_i \in\oQ^*, \ d_1, g \in\Z\setminus\{0\}, \ d_2,
\ldots,d_h\in\Q, \ g^{-2}>d_2>\ldots>d_h>0$.

\nin  We define $k:=h+3$ and, for the $\epsilon>0$ fixed (which we
may take $<1/2k,$ say), $Q=e^{\epsilon}.$ \ We suppose that there
are infinitely many triples $(m,p,q)$ of integers  with
$0<q<Q^{2m+r}$, \ $m\rightarrow +\infty$ \ and
\begin{equation}\label{leq}
\Big|\dfrac{\sqrt{\alpha(2m+r)} + \beta(2m+r)}{\gamma(2m+r)} -
\frac{p}{q}\Big|\leq \frac{1}{q^k} e^{-\epsilon (2m+r)}.
\end{equation}
\nin We shall eventually obtain a contradiction, which will prove
what we want.

\mni We proceed to define the data for an application of the
Subspace Theorem \ref{sub}. We let $S$ be the finite set of places
of $\Q$ containing the infinite one and all the places dividing
the numerators or the denominators of \ $g$ \ and of \ $d_i, \
i=1,\ldots,h$. We define linear forms in $X_o,\ldots,X_h$ as
follows. For $v\neq\infty$ or for $i\neq 0$ we set simply \
$L_{i,v}=X_i$. \ We define the remaining form
$$L_{0,
\infty}:=X_0-b_{1,r}^{1/2}X_1-b_{2,r}X_2-\ldots-b_{h,r}X_h,$$
where $b_{i,r}=b_i b_{1,r}^{1/2}, \ i=2,\ldots,h.$ \ For each $v$,
these linear forms are clearly independent.

\nin Let $d$ be the minimal integer such that $d_id\in\Z$ for
every $i=1,\ldots,h$ (recall Remark 6). For our choice of the set
$S,$ \ $d$ is a S-unit.

\nin Define \ $e_1:=d_1dg^{-2},$ \ $e_i:=dd_i, \ i=2,\ldots,h.$ \
Note that $e_i \in \Z$ for every $i=1,\ldots,h.$

\nin Set the vector
$$\underline{x}=\underline{x}(m,p,q)=(pd^{2m+r}, \ qe_1^md^{m+r}, \
qd_1^me_2^{2m+r},\ \ldots,\ qd_1^me_h^{2m+r}) \in\Z^{h+1}.$$ \nin
We proceed to estimate the double product $\dprod{v \in S}{}\
\dprod{i=0}{h}|L_{i,v}(\underline{x})|_v.$

\nin We have
\begin{equation}\label{prod0}
\dprod{v \in S}{}\
\dprod{i=0}{h}|L_{i,v}(\underline{x})|_v=|L_{0,\infty}(\underline{x})|\cdot\dprod{i=1}{h}\
\dprod{v \in S}{} |L_{i,v}(\underline{x})|_v \cdot \dprod{v\in
S\setminus\{\infty\}}{}|L_{0,v}(\underline{x})|_v.
\end{equation}

\nin By definition $\dprod{v\in
S}{}|L_{1,v}(\underline{x})|_v=\dprod{v\in
S}{}|qe_1^md^{m+r}|_v\leq q$ \ and, for \ $i\geq 2$, \
$\dprod{v\in S}{}|L_{i,v}(\underline{x})|v=\dprod{v\in
S}{}|qd_1^me_i^{2m+r}|_v\leq q$, \ since $d, \ d_1$ and the $e_i$
are S-units for every $i$ (which implies that $\dprod{v\in
S}{}|d|_v=\dprod{v\in S}{}|d_1|_v=\dprod{v\in S}{}|e_i|_v=1$) and
since for the positive integer $q,$ \ $\dprod{v\in S}{}|q|_v\leq
q$ \ holds. This means that
\begin{equation} \label{prod1} \dprod{i=1}{h}\ \dprod{v \in S}{}
|L_{i,v}(\underline{x})|_v\leq q^{h}.
\end{equation}
\nin Moreover,
\begin{equation*} \dprod{v\in
S\setminus\{\infty\}}{}|L_{0,v}(\underline{x})|_v=\dprod{v\in
S\setminus\{\infty\}}{}|pd^{(2m+r)}|_v=
\end{equation*}
\begin{equation}\label{prod2} =\dprod{v\in S\setminus\{\infty\}}{}|p|_v\cdot\dprod{v\in
S\setminus\{\infty\}}{}|d^{(2m+r)}|_v \leq d^{-(2m+r)},
\end{equation}
\nin the last inequality holding since $p$ is an integer and $d$
is a S-unit.

\sni Finally we have

\sni $|L_{0,\infty}(\underline{x})|=
d^{2m+r}\big|p-q\big(b_{1,r}^{1/2}d_1^mg^{-2m}+b_{2,r}d_1^md_2^{2m+r}+\ldots+
b_{h,r}d_1^md_h^{2m+r}\big)\big|=$

\ssk $=qd^{2m+r}\Big|\eta_r(2m+r)-\dfrac{p}{q}\Big|,$

\sni which, combined with (\ref{prod0}), (\ref{prod1}) and
(\ref{prod2}), gives
\begin{equation}
\dprod{v\in S}{}\ \dprod{i=0}{h}|L_{0,v}(\underline{x})|_v\leq
q^{h+1}\Big|\eta_r(2m+r)-\dfrac{p}{q}\Big|.
\end{equation}
\nin Since $q^k<Q^{k(2m+r)}=e^{(2m+r)k\epsilon},$ we have
$q^{-k}e^{-(2m+r)\epsilon}>e^{-(2m+r)(k+1)\epsilon},$ which means
that $q^{-k}e^{-(2m+r)\epsilon}>t^{2m+r}$ (recall that
$\epsilon<1/2k, \ k\geq 3$ \ and \ $t=1/9$). Thus, for a certain
constant $l>0,$ we have

\sni $\Big|\eta_r(2m+r)-\dfrac{p}{q}\Big|\leq
\Big(\Big|\dfrac{p}{q}-\dfrac{\sqrt{\alpha(2m+r)} +
\beta(2m+r)}{\gamma(2m+r)}\Big|+$

\ssk $+\Big|\dfrac{\sqrt{\alpha(2m+r)} +
\beta(2m+r)}{\gamma(2m+r)}-\eta_r(2m+r)\Big|\Big)
\leq\Big(\dfrac{1}{q^k}e^{-(2m+r)\epsilon}+lt^{2m+r} \Big)\leq$

\ssk $\leq\dfrac{2}{q^k}e^{-(2m+r)\epsilon}.$

\sni This means that $\dprod{v\in S}{}\
\dprod{i=0}{h}|L_{0,v}(\underline{x})|_v\leq
2q^{h+1-k}e^{-(2m+r)\epsilon}\leq e^{-(2m+r)\epsilon},$ since we
have $k=h+3$. \ Also, \ $\dmax{0\leq i\leq h}|x_i| \simeq
qe_1^md^{m+r} \leq Q^{2m+r}e_1^md^{m+r}.$

\nin Hence, choosing \ $\delta>0$, \ $\delta < \dfrac {\epsilon}
{\log (Q^2e_1d)},$ \ we get, for $m$ large,
$$\dprod{v\in S}{}\ \dprod{i=0}{h}|L_{0,v}(\underline{x})|_v\leq
e^{-(2m+r)\epsilon}< (Q^{2m+r}e_1^md^{m+r})^{-\delta}\leq
(\dmax{0\leq i\leq h}|x_i|)^{-\delta},$$ \nin i.e. the inequality
of the Subspace Theorem \ref{sub} is verified.

\nin This implies that the vectors
$$\underline{x}=\underline{x}(m,p,q)=(pd^{2m+r}, \ qe_1^md^{m+r}, \
qd_1^me_2^{2m+r},\ \ldots,\ qd_1^me_h^{2m+r}) \in\Z^{h+1}$$ \nin
are contained in a finite set of proper subspaces of $\Q^{h+1}.$
In particular, there exists a fixed subspace, say of equation \
$z_oX_o-z_1X_1-\ldots-z_hX_h=0,$ \ $z_i\in\Q$, containing an
infinity of the vectors in question. We cannot have $z_0=0$, since
this would entail
$z_1e_1^md^{m+r}+z_2d_1^me_2^{2m+r}+\ldots+z_hd_1^me_h^{2m+r}=$

$=d_1^md^{2m+r}(z_1g^{-2m}+z_2d_2^{2m+r}+\ldots+z_hd_h^{2m+r})=0$

\nin for an infinity of $m$; in turn, the fact that $g^{-1}$ and
the $d_i$ are pairwise distinct would imply $z_i=0$ for all $i$, a
contradiction.

\nin Therefore we can suppose that $z_0=1$, and we find that, for
the $m$ corresponding to the vectors in question,
\begin{equation}\label{xi}
\dfrac{p}{q}=d_1^m \Big(z_1g^{-2m}+\dsum{i=2}{h}
z_id_i^{2m+r}\Big)=:\xi(m) \ \in \Q\Sigma_\Q.
\end{equation}
\nin Let us show that actually $\xi \in \Sigma.$ Assume the
contrary; then the minimal positive integer $D$ so that
$D^{m}\xi\in\Sigma$ is $\geq 2.$ But then equation (\ref{xi})
together with Lemma \ref{lem1} implies that $q\gg
2^{m}e^{-m\epsilon}.$ Since this would hold for infinitely many
$m$, we would find $Q\geq
q^{\frac{1}{2m}}\geq\sqrt{2}e^{-\epsilon/2}$, a contradiction
since $Q=e^{\epsilon},$ \ $\epsilon <1/2k$ and $k\geq 3.$

\nin Therefore $\xi \in \Sigma$.

\nin Substituting (\ref{xi}) in (\ref{leq}) we get that there
exists a power sum $\xi\in \Sigma$ such that

$$\Big|\dfrac{\sqrt{\alpha(2m+r)} +
\beta(2m+r)}{\gamma(2m+r)}-\xi(m)\Big|\ll e^{-(2m+r)\epsilon},$$

\nin a contradiction, concluding the proof.

\cvd

\bni {\bf Proof of Corollary \ref{teo1}.} \q We know that
$$l(\alpha-\xi^2)=l((\sqrt{\alpha}-\xi)(\sqrt{\alpha}+\xi))\geq
l(\alpha)^{1/2}$$
\nin holds for every $\xi \in \Sigma$ by
assumption, and that for every $r \in \{0,1\}$
$$|\sqrt{\alpha(2m+r)}+\xi(2m+r)|<2\cdot \max\{\sqrt{\alpha(2m+r)}, \ |\xi(2m+r)|\}.$$
\nin If for a certain $\xi \in \Sigma$ we have \ $|\xi(2m+r)|<k
\cdot \sqrt{\alpha(2m+r)},$ \ for some constant $k>0,$ we get that
for such $\xi \in \Sigma,$
$$|\sqrt{\alpha(2m+r)}-\xi(2m+r)|>\frac{1}{2} \ \min \Big\{1, \frac{1}{k}\Big\}.$$

\nin If for a certain $\xi \in \Sigma$ we have \ $|\xi(2m+r)|\gg
\alpha(2m+r)^{\frac{1}{2}(1+\delta)},$ \ for some $\delta>0,$ we
get
$$|\sqrt{\alpha(2m+r)}-\xi(2m+r)|\gg \alpha(2m+r)^{\frac{1}{2}(1+\delta)}.$$
\nin This proves that there does not exist a power sum \ $\xi \in
\Sigma$ \ and \ $\epsilon>0$ \ such that
$$|\sqrt{\alpha(2m+r)}-\xi(2m+r)|\ll e^{-(2m+r)\epsilon}.$$
\nin Thus we can apply Theorem \ref{teo2} with $\beta=0$ and
$\gamma=1,$ \ and get the conclusion.

\cvd

\bni {\bf Proof of Corollary \ref{cor1}.} \q For notation and
basic facts about continued fractions we refer to \cite{K} and
\cite[Ch. I]{S1}.

\nin Let us suppose by contradiction that there exists an integer
$R>0$ and an infinite set $A\subseteq \N$ such that for $n\in A$
we have
$\sqrt{\alpha(n)}=[a_o(n);\overline{a_1(n),\ldots,a_R(n)}].$

\nin Let $p_i(n)/q_i(n),$ \ $i=0, 1, \ldots,$ \ with $q_0(n)=1$,
be the (infinite) sequence of the convergents of the continued
fraction for $\sqrt{\alpha(n)}$. We recall the relation
$\big|\sqrt{\alpha(n)}-\frac{p_i(n)}{q_i(n)}\big|<(a_{i+1}(n)q_i(n)^2)^{-1}$,
\ for $i\geq 0,$ \ which implies that
\begin{equation}\label{eqcoef}
a_{i+1}(n)<\Big|\sqrt{\alpha(n)}-\frac{p_i(n)}{q_i(n)}\Big|^{-1}q_i(n)^{-2}
\end{equation}
\nin holds for every \ $i\geq 0.$

\nin Since $\alpha$ satisfies the assumptions for Corollary
\ref{teo1}, for some $\epsilon>0$ to be fixed later there exist \
$k>2$ \ and \ $Q=e^{\epsilon}>1$ \ as in the statement.

\nin Define now the increasing sequence \ $c_0, c_1, \ldots$ \ by
\ $c_0=0,$ \ and \ $c_{r+1}=(k+1)c_r+1,$ \ and choose a positive
number \ $\rho <c_R^{-1}\log Q,$ \ so \ $e^{c_R \rho}<Q.$

\nin Proceeding by induction as in the proof of Corollary 1 in
\cite{C-Z}, it can be shown that for every \ $i=0,\ldots,R,$ and
for large $n,$ \ we have \ $q_i(n)<e^{c_i \rho n},$ \ which means
that $q_i(n)<Q^n$ for every $i=0,\ldots,R$ and $n$ large. Thus, we
can apply Corollary \ref{teo1} with $p=p_i(n)$, $q=q_i(n)$, and
$\epsilon>0$ to be chosen later. \ Recalling that
$Q=e^{\epsilon}$, from (\ref{eqcoef}) we get that, for all $n$ but
finitely many, the inequality
\begin{equation}\label{eq1cor2}a_{i+1}(n)< \Big|\sqrt{\alpha(n)}-\frac{p_i(n)
}{q_i(n)}\Big|^{-1}q_i(n)^{-2}\leq q_i(n)^k
e^{n\epsilon}<Q^{kn}e^{n\epsilon}=e^{n(k+1)\epsilon},
\end{equation}
\nin holds for every $i=0,\ldots,R$ \ and \ $\epsilon>0.$

\nin Taking \ $\delta:=(k+1)\epsilon$ \ we can rewrite the above
inequality as
\begin{equation}\label{coef}
a_{i}(n)<e^{n\delta},
\end{equation}
\nin for $i=0,\ldots,R$ and for all $n$ but finitely many.

\nin Let us consider from now on $n\in A$ such that \
$a_{i}(n)<e^{n\delta}$ \ holds.

\nin From well known results of the theory of continued fractions
(see \cite{K}) we get that for every $n$,
\begin{equation}\label{beta}
\sqrt{\alpha(n)}=a_0(n)+\dfrac{1}{\beta(n)},
\end{equation}
\nin where $\beta(n)$ has the continued fraction expansion
$$\beta(n)=[\overline{a_1(n),\ldots,a_R(n)}].$$ \nin This means that $\beta(n)$
satisfies the quadratic equation
$$\beta(n)=[a_1(n),\ldots,a_R(n),\beta(n)],$$ \nin that can be
rewritten as
\begin{equation}\label{bet1}
q'_R(n)\beta(n)^2+(q'_{R-1}(n)-p'_R(n))\beta(n)-p'_{R-1}(n)=0,
\end{equation}
where
$p'_i(n)/q'_i(n)=[a_1(n),\ldots,a_i(n)].$

\nin This means that the integers $p'_{R-1}(n), \ p'_R(n), \
q'_{R-1}(n)$ and $q'_R(n)$ appearing in (\ref{bet1}) are all \
$\ll (\dmax{1\leq i\leq R} a_i(n))^R.$

\nin From (\ref{coef}) it follows that \ $\dmax{1\leq i\leq R}
a_i(n)<e^{n\delta},$ \ which implies that $p'_{R-1}(n), \ p'_R(n),
\ q'_{R-1}(n)$ and $q'_R(n)$ are all \ $\ll e^{Rn\delta}.$

\nin Taking the trace of both terms of (\ref{beta}) we get that
for infinitely many $n$
\begin{equation}\label{somma}
2a_0(n)=\frac{q'_{R-1}(n)-p'_R(n)}{p'_{R-1}(n)}.
\end{equation}
\sni Estimating the height on both sides of (\ref{somma}), on the
left side we get
$$H(2a_0(n))=2a_0(n)=2\lfloor \sqrt{\alpha(n)}
\rfloor \gg 2^{n/2}$$ \nin (since $\alpha$ can be supposed a
non-constant power sum), while on the right side we have
$$H\Big(\dfrac{q'_{R-1}(n)-p'_R(n)}{p'_{R-1}(n)}\Big)\ll
\max{}\{q'_{R-1}(n), p'_R(n), p'_{R-1}(n)\} \ll \ e^{Rn\delta}$$
\nin (since $q'_{R-1}(n), \ p'_R(n),$ \ and \ $p'_{R-1}(n)$ \ are
integers), getting a contradiction choosing $\delta<\dfrac{\ln
2}{2R},$ \ i.e. \ $\epsilon < \dfrac{\ln 2}{2(k+1)R}.$

\cvd

\bni {\bf Proof of the Main Theorem \ref{cor2}.} \q  The case of
$\alpha$ constant is trivial; thus we can suppose $\alpha$ to be
non constant for the rest of the proof.

\nin For $r \in \{0, 1\}$ fixed, let
$$\sqrt{\alpha(2m+r)}=[a_0(m);a_1(m),a_2(m),\ldots]=
[a_0(m);\overline{a_1(m),\ldots,a_{R(m)}(m)}]$$ \nin be the
continued fraction expansion for $\sqrt{\alpha(2m+r)},$ and let \
$p_i(m)/q_i(m),$ \ $i=0, 1, \ldots,$ \ with $q_0(m)=1$, be the
(infinite) sequence of its convergents. If $m \in A,$ \ we have
$R(m)=R.$

\nin We recall that the relations \ $a_R(m)=2a_0(m),$ for every $m
\in A$ (if $R>0$), and
\begin{equation}\label{eq0cor2}
a_{i+1}(m)<\Big|\sqrt{\alpha(2m+r)}-\frac{p_i(m)}{q_i(m)}\Big|^{-1}q_i(m)^{-2},
\end{equation}
\nin for every \ $i\geq 0$ and $m \in \N,$ \ hold.

\nin By our present assumption, the hypothesis of Corollary
\ref{cor1} cannot hold for $\alpha$ and for the fixed $r,$ since
the period of the continued fraction for $\sqrt{\alpha(n)}$ cannot
tend to infinity for $n\longrightarrow +\infty.$ This means that
for a certain $\rho
>0,$ there exists a power sum $\eta \in \Sigma$ such that
\begin{equation}\label{eq1cor2}
|\alpha(2m+r) - \eta(m)^2|\ll \alpha(2m+r)^{1/2-\rho}.
\end{equation}

\nin From (\ref{eq1cor2}) it follows
\begin{equation}\label{eq10cor2}
|\sqrt{\alpha(2m+r)} - \eta(m)|\ll \alpha(2m+r)^{-\rho}<1,
\end{equation}
\nin the last inequality holding for $m \in \N$ large. \nin Since
$\alpha$ has integral coefficients, there exists $\eta$ satisfying
(\ref{eq10cor2}) having the same pro\-perty; this means that
$\eta(m)$ is an integer for every $m$. Since $\eta(m)$ is an
integer and since (\ref{eq10cor2}) holds, it follows that
\begin{equation}\label{eq1bcor2}
a_0(m)=\lfloor\sqrt{\alpha(2m+r)}\rfloor\in \{\eta(m),\
\eta(m)-1\}
\end{equation}
\nin for every $m \in \N$ large enough.

\nin We claim that either \ $a_0(m)=\eta(m)$ \ or \
$a_0(m)=\eta(m)-1$ \ for all \ $m \in \N$ \ large enough. \ In
fact, $a_0(m)=\eta(m)$ when $\alpha(2m+r) - \eta(m)^2\geq0,$ while
$a_0(m)=\eta(m)-1$ \ when \ $\alpha(2m+r) - \eta(m)^2<0,$ and just
one of the above inequalities can hold for all $m$ large, since
$\alpha$ and $\eta$ are power sums. This proves that for $m \in \N
$ large enough $a_0(m)$ is a power sum in $\Z\Sigma_\Z$.

\nin If $R=0,$ the proof is complete.

\sni Note that since $\alpha$ was supposed to be non constant, \
also \ $a_o(m)$ \ is non constant.

\bni Consider from now on $R>0,$ and suppose by contradiction that
there exists $h\in \N,$ \ $1\leq h\leq R,$ \ such that for \ $m
\in A$ \ large enough, \ $a_i(m)$ can be parameterized by a power
sum in $\Z\Sigma_\Z$ \ for \ $i=0,\ldots,h-1,$ \ but not for
$i=h.$

\nin The case \ $h=R$ \ can be excluded, \ since for \ $m \in A$ \
we have \ $a_R(m)=2a_0(m)\in \Z\Sigma_\Z.$

\nin Put $a(m):=[a_0(m);a_1(m), \ldots,
a_{h-1}(m)]=\dfrac{p_{h-1}(m)}{q_{h-1}(m)} \in \Q.$

\nin Since \ $a_i(m)\in \Z\Sigma_\Z$ \ for every \
$i=0,\ldots,h-1,$ \ the relation
\begin{equation}\label{eq2cor2}
|\sqrt{\alpha(2m+r)}-a(m)|^{-1}=\dfrac{\sqrt{\gamma(m)}+\tau(m)}{\xi(m)}=:\alpha_h(m)
\end{equation}
\nin holds for every $m \in A$ large enough, and for certain power
sums \ $\gamma,\tau$  and \ $\xi\in\Z\Sigma_\Z,$ \ $\xi$ not
identically zero.

\mni We claim that for every $\epsilon>0$ there does not exist a
power sum $\zeta\in\Sigma$ such that
\begin{equation}\label{eq11cor2}
\Big| \alpha_h(m)-\zeta(m)\Big|\ll e^{-(2m+r)\epsilon}.
\end{equation}
\nin In fact, if such a power sum would exist, in view of
(\ref{eq11cor2}), we would have $$\Big|
\alpha_h(m)-\zeta(m)\Big|<1$$ \nin for $m \in A$ large enough,
which implies that $$a_h(m)=\lfloor\alpha_h(m)\rfloor\in
\{\lfloor\zeta(m)\rfloor -2, \ \lfloor\zeta(m)\rfloor -1, \
\lfloor\zeta(m)\rfloor\},$$ \nin for $m \in A$ large enough. But
since $\zeta$ has integral roots and rational coefficients, there
exist arithmetic progressions $A_s=\{m=tm'+s, \ m' \in \N\},$ \
for \ $s=0,\ \ldots, \ t-1$ \ and some $t \in \N,$ \ such that
$\lfloor\zeta(m)\rfloor$ can be parameterized by a power sum in
$\Z\Sigma_\Z$ for all $m\in A$ in any of such progressions. Choose
a progression, say $A_1,$ that contains infinitely many elements
$m \in A.$ Let us notice that the set \ $A$ \ in the statement of
the present Theorem can be substituted without losing generality
by any of its infinite subsets ($A$ is just an infinite set for
which $R(m)=R,$ \ and not the set of all $m$ for which $R(m)=R$).
\ Substituting the set $A$ in the statement of the present Theorem
by the (still infinite) set $A\cap A_1,$ which for simplicity of
notation we will call $A$ again, we would get that for all $m \in
A$ large enough, $a_h(m)$ can be parameterized by a power sum in
$\Z\Sigma_\Z,$ a contradiction proving that $\alpha_h$ satisfies
the assumption of Theorem \ref{teo2}.

\mni By the definition of $\alpha_h(m),$ the length of the period
of its continued fraction is $R$ again. Let
$$\alpha_h(m)=[a_0'(m);\overline{a_1'(m),\ldots,a_R'(m)}],$$
and let \ $p_i'(m)/q_i'(m),$ \ $i=0, 1, \ldots,$ \ with
$q_0'(m)=1$, be the (infinite) sequence of its convergents.

\nin We have the relations \ $a_i'(m)=a_{i+h}(m)$ \ for \ $i+h\leq
R$, \ \ $a_i'(m)=a_{i+h-R}(m)$ \ for \ $i+h>R$, \ and
\begin{equation}\label{eq2'cor2}
a_{i+1}'(m)<\Big|\alpha_h(m)-\dfrac{p_i'(m)}{q_i'(m)}\Big|^{-1}
\end{equation}
\nin for every $i\geq 0.$

\nin Since $\alpha_h$ satisfies the assumption for Theorem
\ref{teo2}, for some $\epsilon>0$ to be fixed later there exist \
$k \geq 3$ \ and \ $Q=e^\epsilon>1$ \ as in that statement.

\nin As in the proof of Corollary \ref{cor1}, we have again the
inequality \ $q_i'(m)<Q^{2m+r},$ \ which holds for every \
$i=0,\ldots,R$ and $m$ large, \ i.e. we can apply Theorem
\ref{teo2} to $\alpha_h(m)$ with $p=p_i'(m),$ \ $q=q_i'(m)$ \ and
some \ $\epsilon>0$ to be fixed later. We get that for every
$i\geq 0$ and for $m \in A$ large enough,
\begin{equation}\label{eq3cor2}
\Big|\alpha_h(m)-\dfrac{p_i'(m)}{q_i'(m)}\Big|\geq q_i'(m)^{-k}
e^{-(2m+r)\epsilon}.
\end{equation}
\nin Recalling that \ $0<q_i'(m)<Q^{2m+r}=e^{(2m+r)\epsilon}$, \
for every \ $i=0, \ldots, R,$ and considering the inequality
(\ref{eq3cor2}) for $i=R-h-1,$ together with (\ref{eq2'cor2}), we
have
\begin{equation*}
a_{R}(m)=a_{R-h}'(m)\leq
\Big|\alpha_h(m)-\dfrac{p_{R-h-1}'(m)}{q_{R-h-1}'(m)}\Big|^{-1}\leq
q_{R-h-1}'(m)^{k} \ e^{(2m+r)\epsilon}<
\end{equation*}
\begin{equation}\label{eq4cor2}
<Q^{(2m+r)k}e^{(2m+r)\epsilon}=e^{(2m+r)(k+1)\epsilon}=e^{(2m+r)\epsilon'},
\end{equation}
\nin for $\epsilon'=(k+1)\epsilon.$

\nin Choosing $\epsilon<\frac{\ln2}{2(k+1)}$ \ (i.e.
$\epsilon'<\frac{\ln2}{2}$), \ we get that $$a_R(m)\ll
2^{m(1-\delta)},$$
\nin for some $\delta>0.$

\nin Recalling that $a_0(m)\in\Z\Sigma_\Z$ is non constant, from
the relation $$a_R(m)=2a_0(m)\gg 2^m$$ we get a contradiction,
proving that the relation (\ref{eqcor2}) holds for every $m \in
A,$ except finitely many.

\bni It remains to show that (\ref{eqcor2}) holds for every $m \in
\N,$ except finitely many. We will proceed by contradiction.

\nin We have already proved that $a_0(m)=\beta_0(m)$ for every $m
\in \N$ large enough.

\nin Suppose that for some $u>0,$ \ $a_i(m)=\beta_i(m)$ for every
$i=0, \ldots, u-1$ \ and for every $m \in \N$ except finitely
many, but $a_u(m) \neq \beta_u(m)$ for infinitely many $m \in \N$
\ (we define \ $\beta_{aR+b}(m):=\beta_b(m),$ for $a\in \N$ and
$0\leq b<R$).

\nin Let $a'(m):=[\beta_0(m), \ldots, \beta_{u-1}(m)].$

\nin We know that for $m \in \N$ large enough,
$$|\sqrt{\alpha(2m+r)}-a'(m)|^{-1}=\dfrac{\sqrt{\gamma'(m)}+\tau'(m)}{\xi'(m)},$$
\nin for certain $\gamma', \eta', \xi' \in \Z\Sigma_\Z, \q \xi'$
not identically zero.

\nin For $m \in A$ large enough we have
\begin{equation}\label{eq5cor2}
\beta_u(m)=a_u(m)=\lfloor | \sqrt{\alpha(2m+r)}-a'(m)|^{-1}
\rfloor=\Big\lfloor
\dfrac{\sqrt{\gamma'(m)}+\tau'(m)}{\xi'(m)}\Big\rfloor,
\end{equation}
\nin which means that both the inequalities
\begin{equation}\label{eq6cor2}
\dfrac{\sqrt{\gamma'(m)}+\tau'(m)}{\xi'(m)} - \beta_u(m)\geq 0
\end{equation}
\nin and
\begin{equation}\label{eq7cor2}
\dfrac{\sqrt{\gamma'(m)}+\tau'(m)}{\xi'(m)} - \beta_u(m)<1
\end{equation}
\nin hold for $m \in A$ large enough.

\nin The inequalities (\ref{eq6cor2}) and (\ref{eq7cor2}) can be
rewritten as
\begin{equation}\label{eq8cor2}
\gamma'(m)-(\beta_u(m)\xi'(m)-\tau'(m))^2\geq 0
\end{equation}
and
\begin{equation}\label{eq9cor2}
\gamma'(m)-(\xi'(m)+\xi'(m)\beta_u(m)-\tau'(m))^2<0
\end{equation}
\nin respectively.

\nin Since $\beta_u, \gamma', \tau', \xi' \in \Z\Sigma_\Z$ are
power sums, both the inequalities (\ref{eq8cor2}) and
(\ref{eq9cor2}) can hold either for every $m \in \N$ except
finitely many, or just for a finite set of $m.$ Since we know that
they hold for an infinite subset of $A,$ they must hold for every
$m \in \N,$ except at most finitely many, \ i.e.
$\beta_u(m)=a_u(m)$ for every $m \in \N$ except finitely many, a
contradiction proving that
$$\sqrt{\alpha(2m+r)}=[\beta_0(m);\overline{\beta_1(m),\ldots,\beta_R(m)}]$$
\nin for every \ $m\in \N,$ \ apart from finitely many exceptions.

\cvd

\bni

\nin \textbf{Acknowledgement:} \q The author is grateful to Prof.
U. Zannier for bringing this problem to his attention and to Prof.
P. Corvaja for valuable comments and remarks.

\bni

\bsk

\bsk

\bsk

\noindent{\sc Amedeo Scremin\\
Institut f\"ur Mathematik\\
TU Graz\\
Steyrergasse 30\\
A-8010 Graz, Austria\\
e-mail: {\sf scremin@finanz.math.tu-graz.ac.at}}

\end{document}